\ifcase1
\or
\documentclass[10pt,a4paper,twocolumn]{article}
\or
\documentclass[12pt,a5paper]{article}
\fi

\usepackage{euler,amsmath,defns}
\IfFileExists{ajr.sty}{\usepackage{ajr}}{}
\allowdisplaybreaks

\newcommand{\ee}{{\textsc{e}}}
\newcommand{\eee}{{\bar\ee}}
\newcommand{\uu}{{\bar u}}

\newcommand{\nuu}{{\bar \nu}}
\newcommand{\rnuu}{R_\nuu}
\newcommand{\ros}{{\dot\varepsilon}}
\newcommand{\vel}{{\vec u}}
\renewcommand{\gr}{}\renewcommand{\hat}{} 

\title{The inertial dynamics of thin film flow of non-Newtonian
fluids}

\author{A.~J. Roberts\thanks{Computational Engineering and Science
Research Centre, Department of Mathematics \& Computing, University of
Southern Queensland, Toowoomba, Queensland 4352, Australia.
\protect\url{mailto:aroberts@usq.edu.au}}}

\begin{document}
    
\maketitle
    
\begin{abstract}
Consider the flow of a thin layer of non-Newtonian fluid over a solid
surface.  I model the case of a viscosity that depends nonlinearly on
the shear-rate; power law fluids are an important example, but the
analysis here is for general nonlinear dependence.  The modelling
allows for large changes in film thickness provided the changes occur
over a large enough lateral length scale.  Modifying the surface
boundary condition for tangential stress forms an accessible base for
the analysis where flow with constant shear is a neutral critical mode,
in addition to a mode representing conservation of fluid.
Perturbatively removing the modification then constructs a model for
the coupled dynamics of the fluid depth and the lateral momentum.  For
example, the results model the dynamics of gravity
currents of non-Newtonian fluids even when the flow is not very slow.
\end{abstract}




\section{Introduction}

Consider the two~dimensional flow of a thin layer of fluid over a flat
substrate.  The fluid of thickness~$\eta(x,t)$ spreads with mean
lateral velocity~$\uu(x,t)$.  Suppose the fluid has the non-Newtonian,
power law, stress-strain relation that the $\text{stress} \propto
(\text{strain-rate})^s$ for some fixed exponent~$s$: the exponent $s=1$
for a Newtonian fluid; $s<1$ is shear thinning; and $s>1$ is shear
thickening.  Such a power law is sometimes called Ostwald's or Norton's
constitutive relation~\cite{Gratton99}.  Then the systematic analysis
developed in this article supports the nondimensional model
\begin{eqnarray} &&
    \D t\eta +\D x{}\left[ \eta\uu \right] =0\,,
    \label{eq:ghup0}
    \\&&
    \re\left[ \D t\uu +\frac{167-25/s}{96}\uu\D x\uu 
    +\frac{25-13/s}{96\eta}\uu^2\D x\eta \right]
    \nonumber\\&&\approx
    - \frac {5(25-1/s)c_s}{48\sqrt2\eta}\left(\frac{\sqrt2\uu}{\eta}\right)^s
    \nonumber\\&&\quad{}
    + \frac{19+1/s}{24}\left(g_1-g_2\D x\eta \right) \,,
    \label{eq:gup0}
\end{eqnarray}
where $\re$~is the nondimensional Reynolds number, $c_s$~is the
coefficient of proportionality in the nonlinear stress-strain relation,
and where $g_1$~and~$g_2$ are the nondimensional components of gravity
along and normal to the flat substrate.  Fluid is conserved
through~\eqref{eq:ghup0}.  The momentum equation~\eqref{eq:gup0}
incorporates effects of inertia, self-advection, bed drag and
gravitational forcing; the dependence of the coefficients upon~$s$
models the subtle effects of the power law rheology.

This model not only applies to the flow of simple liquids, it applies
to: gravity currents of suspensions with medium to high volume
fractions as these are non-Newtonian~\cite{Stickel05}; power law
rheologies are used to model ice flow~\cite[e.g.]{Payne99,
Wilchinsky04} and at even a few metres per year the Reynolds number is
significant for a thick glacier; and a modified model would apply to
turbulent flow as the Smagorinsky large eddy closure of turbulence
corresponds to the shear thickening case of exponent $s=2$
\cite[Eqn.~(6), e.g.]{Ozgokmen07}.  This article puts models such as
\eqref{eq:ghup0}--\eqref{eq:gup0} within the sound support of modern
dynamical systems theory, Section~\ref{sec:smwcmt}, to empower us to
systematically control error, assess domains of validity, and to
systematically account for further physical effects.

The analysis here encompasses not only power law fluids but a general
nonlinear dependence of the stress upon the strain-rate as codified in
Section~\ref{sec:demf}.  In contrast, almost all previous thin fluid
film modelling use only a power law dependence.  Some industrial
plastics have a complicated non-monotonic dependence~\cite{Bird95b}
that cannot be represented by a simple power law.  Similarly, dense
suspensions often have non-monotonic dependence~\cite{Stickel05}.  The
resultant model derived in Section~\ref{sec:lomd} also applies to such
complicated industrial plastics and dense suspensions.

The lubrication approximation of very slow flow, low Reynolds number,
underpins previous theoretical models for non-Newtonian thin fluid
films: Perazzo \& Gratton~\cite{Perazzo03} and Betelu \&
Fontelos~\cite{Betelu04} examined flow with surface tension; this
followed experiments compared with travelling waves and similarity
solutions by Gratton, Minotti \& Mahajan~\cite{Gratton99}.  Gratton et
al.\ comment ``the differences between Newtonian and non-Newtonian
currents are significant and can clearly be observed in experiments''.
But the lubrication approximation, that creates models expressed only
in terms of the fluid thickness~$\eta(x,t)$, does not model inertia and
so cannot resolve any wave-like dynamics.  To model faster flows,
potentially with wave effects, we must resolve the dynamics of both the
fluid thickness and a measure of horizontal momentum~\cite{Roberts94c,
Roberts99b}, we used $\eta$~and~$\uu$
in~\eqref{eq:ghup0}--\eqref{eq:gup0}.  For example, Harris et
al.~\cite{Harris01} modelled particle driven gravity currents using
shallow water equations that resolve the dynamics of both the fluid
thickness and the mean lateral velocity.  However, such modelling of
essentially dissipative flows, albeit dissipative via turbulence, by
the laminar inviscid foundation of shallow water equations appears a
contradiction that demands resolution.  This article shows how such
models of non-Newtonian fluid flow may be put on a sound mathematical
basis to empower  accurate physical forecasts.

\section{Differential equations to model non-Newtonian flow}
\label{sec:demf}

Let the incompressible fluid have thickness~$\eta(x,t)$, constant
density~$\rho$, a nonlinear rheology, and let the fluid flow with some
varying velocity field~$\vel =(u,v)=(u_1,u_2)$ and pressure field~$p$.

\paragraph{Nonlinear constitutive relation}
Define the strain-rate tensor \cite{Gratton99,
Stickel05}\footnote{Some, such as Betelu \& Fontelos~\cite{Betelu04},
use double this tensor.}
\begin{equation}
	\ros_{ij} =\rat12\left( \D{x_j}{u_i} +\D{x_i}{u_j} \right) \,,
\end{equation}
where $x_1=x$ and $x_2=y$ are distances along and normal to the solid
substrate, respectively.  Then the stress tensor for the
fluid is $\sigma_{ij} =-p\delta_{ij} +2\rho\nu \ros_{ij}$\,: for a
Newtonian fluid the kinematic viscosity~$\nu$ is constant; but when the
kinematic viscosity varies with strain-rate then we model shear
thickening or shear thinning non-Newtonian fluids.

The important class of non-Newtonian fluids that we address has
viscosity which depends only upon the magnitude~$\ros$ of the second
invariant of the strain-rate tensor~\cite{Betelu04}:
\begin{equation} 
    \ros^2=\sum_{i,j}\ros_{ij}^2\,.
    \label{eq:cons}
\end{equation}
For example, Bird et al.~\cite[see~\cite{Betelu04}]{Bird77} report that
a solution of 0.5\%~Hydroxyethylcellulose is shear thinning: at~$20\C$
the solution has viscosity $\mu=m{\ros}^{s-1}$ for exponent $s=1/1.96$
and coefficient $m=0.84\N\secs^{s}/\m^2$\,.

\paragraph{Partial differential equations} 

Make equations nondimensional with respect to some velocity scale, a
typical fluid thickness, and the fluid density.  The nondimensional
\pde{}s for the incompressible, two~dimensional, fluid flow are firstly
the continuity equation
\begin{equation}
    \divv\vel =\D xu+\D yv=0\,,
\end{equation}
and secondly the momentum equation
\begin{equation}
    \re\left( \D t{\vel } +\vel \cdot\grad\vel  \right)
    =-\grad p +\divv\vec\tau +\gr\vec{\hat g}\,,
\end{equation}
where $\re$~is the appropriate Reynolds number, $\vec\tau$~is the
nondimensional deviatoric stress tensor, and $\vec g=(g_1,g_2)$ is the
nondimensional forcing of gravity.  For a fluid with a nonlinear
stress-strain relation, the nondimensional deviatoric stress tensor
\begin{equation}
    \tau_{ij}=2\nu(\ros)\ros_{ij}
    =\nu(\ros)\left( \D{x_j}{u_i} +\D{x_i}{u_j} \right)\,.
\end{equation}

\paragraph{Boundary conditions}
Solve these \pde{}s  with nondimensional boundary conditions:
\begin{itemize}
	\item on the bed of no-slip,\footnote{If modelling turbulent flows
	by a large eddy closure, we may justifiably replace this no-slip
	bed condition by a mixed boundary condition on the lateral
	velocity: $u\propto \D yu$\,.}
    \begin{equation}
        \vel =\vec 0\quad\text{on}\quad y=0\,;
        \label{eq:noslip}
    \end{equation}

    \item  the kinematic condition on the free-surface of 
    \begin{equation}
        \D t\eta+u\D x\eta=v\quad\text{on}\quad y=\eta\,;
    \end{equation}

    \item  the stress normal to the free surface comes from constant 
    environmental pressure and surface tension, that is,
    \begin{eqnarray}&&
		-p+\frac1{1+\eta_x^2}\left(\tau_{22} -2\eta_x\tau_{12}
		+\eta_x^2\tau_{11} \right)
		\nonumber\\&&{}
		= \frac{\we\eta_{xx}}{(1+\eta_x^2)^{3/2}}
		 \quad\text{on}\quad y=\eta\,,
    \end{eqnarray}
    where $\we$~is a nondimensional Weber number characterising the
    importance of surface tension;
    
    \item and there must be no tangential stress at the free surface,
    \begin{equation}
		(1-\eta_x^2)\tau_{12}+\eta_x(\tau_{22}-\tau_{11})=0
		\quad\text{on}\quad y=\eta\,.
        \label{bc:tt}
    \end{equation}
This boundary condition of zero tangential stress implicitly is
effectively one of zero shear at the surface; this is not appropriate
for material with a finite yield stress.  Here we assume the fluid
yields for arbitrarily small stress.
\end{itemize}

\section{Centre manifold theory supports the modelling}
\label{sec:smwcmt}

This section shows how to place models such as
\eqref{eq:ghup0}--\eqref{eq:gup0} on a sound theoretical base.
Artificially modify the zero tangential stress free surface
condition~\eqref{bc:tt} to have an artificial forcing proportional to
the local velocity, a forcing which we later remove by evaluating at
parameter $\gamma=1$\,:
\begin{eqnarray}&&
    (1-\rat16\gamma)\left[
    (1-\eta_x^2)\tau_{12}+\eta_x(\tau_{22}-\tau_{11})
    \right] 
    \nonumber\\&&{}
    = (1-\gamma)\frac{\nu(\ee)}{\eta} u
    \quad\text{on}\quad y=\eta\,.
    \label{bc:tta}
\end{eqnarray}
Evaluated at $\gamma=1$ this artificial right-hand side becomes zero so
the boundary condition~\eqref{bc:tta} reduces to the physical boundary
condition~\eqref{bc:tt}.  However, when the parameter $\gamma=0$ and
the lateral gravity and lateral derivatives negligible, $\hat
g_1=\partial_x=0$\,, a neutral mode of the dynamics is the lateral
shear $u=\sqrt2\ee y$ where I define $\ee$~to be proportional to the mean
lateral strain-rate:
\begin{displaymath}
    \ee=\frac1{\sqrt2\eta} \int_0^\eta \D yu\,dy
    =\frac1{\sqrt2\eta} u|_{y=\eta}\,.
\end{displaymath}
This neutral lateral shear mode arises because in pure shear flow
$\tau_{12}=\nu u_y$ and hence the artificial free surface
condition~\eqref{bc:tta} reduces to $\nu u_y =\nu u/\eta$ on
$y=\eta$\,.  Conservation of fluid provides a second neutral mode in
the dynamics.  That is, when $\gamma=\hat g_1=\partial_x=0$ then a two
parameter family of equilibria exists corresponding to some uniform
lateral shear flow, $u=\ee y$\,, on a fluid of any constant
thickness~$\eta$.  For large enough lateral length scales, these
equilibria occur independently at each
location~$x$~\cite[e.g.]{Roberts88a, Roberts96a} and hence the space of
equilibria are in effect parametrised by $\ee(x)$~and~$\eta(x)$.
Provided we can treat lateral derivatives~$\partial_x$ as a modifying
influence, that is provided solutions vary slowly enough in~$x$, centre
manifold theory~\cite[e.g.]{Carr81, Kuznetsov95, Roberts03b} assures us
three things: this space of equilibria is perturbed to a slow manifold,
on which the evolution is slow, that \emph{exists for a finite range}
of $\gamma$~and~$\hat g_1$, and which may be parametrised by the mean
lateral shear~$\ee(x,t)$ and the local thickness of the
fluid~$\eta(x,t)$; the slow manifold is \emph{attractive for all nearby
initial conditions}; and that a formal power series in the
parameters~$\gamma$, $\hat g_1$ and derivatives~$\partial_x$
\emph{approximates} the slow manifold.  That is, the theory supports
the existence, relevance and construction of slow manifold models such
as \eqref{eq:ghup0}--\eqref{eq:gup0}.

\section{Low order models of the dynamics}
\label{sec:lomd}

Computer algebra readily constructs slow manifold models as asymptotic
solutions of the governing differential equations and boundary
conditions.

\subsection{Power law fluids}
For simplicity, suppose the rheology is a nondimensional power law
for the kinematic viscosity, $\nu=c_s\ros^{s-1}$\,.

Computer algebra~\cite[\S3]{Roberts07a} derives that for such a power law
fluid, the evolution of the fluid thickness~$\eta$ and the stress
parameter~$\ee$ is
\begin{eqnarray}
    \D t\eta &=&
    -\D x{} \left[ \left(1+\rat{5\sqrt2}{48 s}\gamma \right)
    \rat12\eta^2\ee \right]
    \nonumber\\&&{}
    +\Ord{\partial_x^2+\hat g_1^2+\gamma^3}\,,
    \label{eq:ghp} \\
    \re\D t\ee &=&
    -\rat52\left(\gamma +\rat{1-1/s}4\gamma^2 \right) c_s\frac{\ee^s}{\eta^2}
    \nonumber\\&&{}
    -\re\sqrt2\left[ \left(\rat38 +\rat{1-8/s}{96}\gamma \right) \eta\ee\D x\ee 
        -\rat1{6s}\gamma\ee^2\D x\eta \right]
    \nonumber\\&&{}
    +\gr \sqrt2\left[ \rat34 -\rat{1+1/s}{16}\gamma \right]
    \eta^{-1}\left(\hat g_1-\hat g_2\D x\eta \right)
    \nonumber\\&&{}
    +\Ord{\partial_x^2+\hat g_1^2+\gamma^3}\,.
    \label{eq:gep}
\end{eqnarray}
The nonlinear rheology primarily appears as a nonlinear drag on the
bed.  However, changes to the vertical profiles of velocity and
pressure due to different power laws affect the coefficients of this
model through their dependence upon exponent~$s$.  

In modelling the flow of thin fluid layers, researchers generally
prefer to use the mean lateral velocity or the lateral fluid flux
instead of the shear parameter~$\ee$.  Using the velocity fields
computed simultaneously with \eqref{eq:ghp}--\eqref{eq:gep} the
computer algebra~\cite[\S3]{Roberts07a} also derives the mean lateral
velocity
\begin{eqnarray*}
   \uu&=&\frac1\eta\int_0^\eta u\,dy 
   = \frac1{\sqrt2}\left( 1+\rat5{24s}\gamma
   +\rat{5(4-1/s)}{288s}\gamma^2 \right)\eta\ee 
   \\&&{}
   +\frac{\eta^2\ee^{1-s}}{sc_s}\left[ \rat\re{160}\ee\D x\ee
   +\gr\rat1{48}\left(\hat g_1-\hat g_2\D x\eta \right) \right]
   +\cdots\,.  
\end{eqnarray*}
Reverting this series to express~$\ee$ in terms of~$\uu$, and
substituting into the model \eqref{eq:ghp}--\eqref{eq:gep} leads to a
model for the coupled evolution of $\eta(x,t)$~and~$\uu(x,t)$.
Evaluating at the physically relevant $\gamma=1$, to remove the
artifice in the surface boundary condition~\eqref{bc:tta}, then gives
the model \eqref{eq:ghup0}--\eqref{eq:gup0} discussed in the
Introduction of this article.

Computer algebra experiments~\cite[\S1]{Roberts07a} suggest that the
convergence of the asymptotic series in~$\gamma$ is markedly improved
by the factor $(1-\rat16\gamma)$ on the left-hand side of the
tangential stress boundary condition~\eqref{bc:tta}.  This factor is
equivalent to an Euler transformation of the asymptotic series.  As
shown in other similar applications~\cite[e.g.]{Roberts94c}, evaluation
at $\gamma=1$ is valid provided the lateral derivatives are small
enough.

Computer algebra~\cite[\S3]{Roberts07a} may construct terms in the
formal power series solutions to higher order in the parameters
$\gamma$, $\hat g_1$ and~$\partial_x$ to generate many valid
approximations of varying orders of accuracy.  For example, to resolve
any effects of surface tension we need to compute terms
in~$\partial_x^2$ that are neglected in \eqref{eq:gup0}
and~\eqref{eq:gep}.  With the support of centre manifold theory,
researchers may choose an approximate model that suits the parameter
regime of their application.

\subsection{More general non-Newtonian fluids} 

We now return to the more general case where the viscosity~$\nu$ of the
fluid depends quite generally upon the magnitude of the
shear-rate~$\ros$, instead of being a simple power law.  In this more
general case the expressions for the modelling are much more
complicated.  The reason is the general nonlinear dependence of
viscosity on strain-rate: for conciseness define
\begin{equation}
    \eee=\frac{\sqrt2\uu}\eta\,,\quad
    \nuu=\nu(\eee)
    \qtq{and}
    \rnuu=\frac1{\nuu+\eee\nuu'}\,,
\end{equation}
where primes on~$\nuu$ denote the derivatives~$d/d\ee$ of the
viscosity~$\nu(\ee)$ and evaluated at $\ee=\sqrt2\uu/\eta$\,.  

The procedure is as for the power law case: computer
algebra~\cite[\S3]{Roberts07a} constructs the slow manifold and
evolution thereon to some order of error; then revert the asymptotic
series to find stress parameter~$\ee$ as a function of mean
velocity~$\uu$; and substitute to express the model in terms of
$\eta$~and~$\uu$.  Conservation of fluid again derives~\eqref{eq:ghup0}
(to any order of error).  The dynamics of momentum then leads to
\newcommand{\fudge}{\hspace{-1.6em}}
\begin{eqnarray}&&\fudge
    \re\D t\uu=
    -\left[\rat{5\gamma}2 +\rat{5\gamma^2}{48}\eee\nuu\rnuu^2 \big(
    2\nuu'+\eee\nuu'' \big) \right]\frac{\nuu\uu}{\eta^2}
    \nonumber\\&&\fudge{}
    -\re\left[\rat74 -\rat{13\gamma}{48} +\rat\gamma{96}\eee\rnuu^2
    \big( 38\nuu\nuu'+12\eee\text{$\nuu'$}^2+13\eee\nuu\nuu'' \big) 
    \right]\uu\D x\uu
    \nonumber\\&&\fudge{}
    -\re\sqrt2\left[\rat18 -\rat\gamma{16}
    +\rat{13\gamma}{192}\eee^2\rnuu^2 \big( 2\text{$\nuu'$}^2-\nuu\nuu''
    \big ) \right]\eee\uu\D x\eta
    \nonumber\\&&\fudge{}
    +\gr\left[\rat34+\rat\gamma{12}
    -\rat\gamma{24}\eee\nuu\rnuu^2 \big( 2\nuu'+\eee\nuu'' \big) \right]
    \left(\hat g_1-\hat g_2\D xh\right)
    \nonumber\\&&\fudge{}
    +\Ord{\partial_x^2+\hat g_1^2+\gamma^3}\,.
    \label{eq:gu}
\end{eqnarray}
Evaluate this equation at $\gamma=1$ to recover a physically relevant
model of the dynamics of lateral momentum.  

The power law model~\eqref{eq:gup0} is just one specific subclass of
the general model~\eqref{eq:gu}: obtain~\eqref{eq:gup0} by the
specific choice of a power law viscosity,
$\nu(\ros)=c_s\ros^{s-1}$\,.

\section{Conclusion}

Following similar modelling for Newtonian thin films~\cite{Roberts94c},
this approach places the modelling of an important class of
non-Newtonian fluids upon the sound basis of centre manifold
theory~\cite[e.g.]{Carr81, Kuznetsov95, Roberts03b}.  This modern
dynamical system foundation empowers us to systematically derive the
novel and accurate models~\eqref{eq:gup0}, \eqref{eq:gep}
and~\eqref{eq:gu} for the lateral momentum of fluids with nonlinear
rheology.

These models of thin fluid flow can be directly applied to flows as
diverse as those of industrial plastics~\cite[e.g.]{Bird95b},
ice~\cite[e.g.]{Payne99, Wilchinsky04}, and medium to dense
suspensions~\cite[e.g.]{Stickel05}.  When you desire more accuracy than
that presented here, computer algebra readily computes higher order
approximations~\cite[\S3]{Roberts07a}.  Modifying the no-slip boundary
condition on the bed,~\eqref{eq:noslip}, will empower the modelling of
turbulent layers of flow over a substrate via large eddy closures.
There are enormous applications for this approach to modelling the
dynamics of laterally extensive layers of fluids.

\bibliographystyle{plain}
\bibliography{new,ajr,bib}

\end{document}